\documentclass[english]{article}
\usepackage[cp850]{inputenc}
\usepackage{babel}
\usepackage{amsmath}
\usepackage{amsfonts,amssymb}
\usepackage[all]{xypic}
\RequirePackage{theorem}%
\RequirePackage{thc}%
\theoremstyle{change} \theoremheaderfont{\normalfont\bfseries}

\def\definitionName{Definition.}
\def\remarkName{Remark.}
\def\remarksName{Remarks.}
\def\noteName{Note.}
\def\notesName{Notes.}
\def\questionName{Question.}
\def\problemName{Problem.}
\def\exampleName{Example.}
\def\examplesName{Examples.}
\def\exerciseName{Exercise.}
\def\exercisesName{Exercises.}
\def\lpName{}
\def\hunchName{Hunch Hunch.}
\def\theoremName{Theorem.}
\def\propositionName{Proposition.}
\def\lemmaName{Lemma.}
\def\corollaryName{Corollary.}
\def\lppName{}
\def\proofName{Proof.}
\def\solutionName{Solution.}

\newcounter{mysubsection}[section]
\newcounter{mysubsubsection}[mysubsection]

\renewcommand{\thesection}{\arabic{section}}%
\renewcommand{\themysubsubsection}{(\thesection.\arabic{mysubsection}.\arabic{mysubsubsection})}%

\theorembodyfont{\upshape}

\newtheorem{remark}[mysubsection]{\remarkName}

\newtheorem{example}[mysubsection]{\exampleName}

\theorembodyfont{\slshape}
\newtheorem{theorem}[mysubsection]{\theoremName}
\newtheorem{proposition}[mysubsection]{\propositionName}
\newtheorem{lemma}[mysubsection]{\lemmaName}

  {\par\addvspace{\smallskipamount}%
    \refstepcounter{mysubsubsection}%
    {\upshape\bfseries%
     \themysubsubsection%
     \hskip\labelsep\lppName}}%
  {\par\addvspace{\smallskipamount}}%

\def\qed{{\unskip\nobreak\hfil\penalty50%
  \hskip2em\hbox{}\nobreak\hfil$\square$%
  \parfillskip=0pt\finalhyphendemerits=0\par}}

\newenvironment{proof}%
  {\par\addvspace{\medskipamount}%
    \upshape%
    {\scshape\bfseries%
    \proofName\hskip\labelsep}}%
  {\qed%
    \addvspace{\medskipamount}}%
  {\par\addvspace{\medskipamount}%
    \upshape%
    {\scshape\bfseries%
    \proofName\hskip\labelsep}#1}%
  {\qed%
    \addvspace{\medskipamount}}%
  {\par\addvspace{\medskipamount}%
    \upshape%
    {\scshape\bfseries%
    \solutionName\hskip\labelsep}}%
  {\qed%
    \addvspace{\medskipamount}}%

\newcommand{\Alg}   {{\rm Alg}}
\newcommand{\End}   {{\rm End}}
\newcommand{\Hom}   {{\rm Hom}}
\newcommand{\Ker}   {{\rm Ker}}
\newcommand{\balf}   {\renewcommand{\theenumi}{(\alph{enumi})}
                      \renewcommand{\labelenumi}{\theenumi}
                      \begin{enumerate}}
\newcommand{\ealf}   {\end{enumerate}
                      \renewcommand{\theenumi}{\arabic{enumi}}
                      \renewcommand{\labelenumi}{\theenumi.}}
\newcommand{\bara}   {\renewcommand{\theenumi}{(\arabic{enumi})}
                      \renewcommand{\labelenumi}{\theenumi}
                      \begin{enumerate} }
\newcommand{\eara}   {\end{enumerate}
                      \renewcommand{\theenumi}{\arabic{enumi}}
                      \renewcommand{\labelenumi}{\theenumi.}}
\newcommand{\brom}   {\renewcommand{\theenumi}{(\roman{enumi})}
                      \renewcommand{\labelenumi}{\theenumi}
                      \begin{enumerate} }
\newcommand{\erom}   {\end{enumerate}
                      \renewcommand{\theenumi}{\arabic{enumi}}
                      \renewcommand{\labelenumi}{\theenumi.}}

\begin{document}
\title{Lie Bracket of Vector Fields in Noncommutative Geometry}
\author{\begin{tabular}{rcl}
 P. Jara\thanks{Research partially supported by DGES BMF2001-2823 and
                FQM--266 (Junta de Andaluc{\'\i}a Research Group).\;
                \texttt{pjara@ugr.es};\;
                \texttt{http://www.ugr.es/local/pjara}}
 &&
 D. Llena\\
 Department of Algebra
 &&
 Department of Geometry\\
 University of Granada
 &&
 University of Almer{\'\i}a\\
 18071--Granada. Spain
 &&
 04120--Almer{\'\i}a. Spain
\end{tabular}}

\date{}
 \maketitle
 \thispagestyle{empty}

\begin{abstract}
The aim of this paper is to avoid some difficulties, related with the Lie bracket, in
the definition of vector fields in a non commutative setting, as they were defined by
Woronowicz, Schm{\"u}dgen--Sch{\"u}ler and Aschieri--Schupp. We extend the definition of vector
fields to consider them as ``\emph{derivations}" of the algebra, through Cartan pairs
introduced by Borowiec. Then, using translations, we introduce the invariant vector
fields. Finally, the definition of Lie bracket realized by Dubois--Violette,
considering elements in the center of the algebra, is also extended to these invariant
vector fields.

\noindent
\underline{Classification numbers}. 2003--\textbf{PACS:} 02.40.Ma, 03.65.Fd%

\end{abstract}

\section*{Introduction}

In the present paper we work on a Hopf algebra $H$ in order to define on it vector
fields as ``\emph{dual}" of differential forms. Recall that differential forms was the
usual way to study non commutative geometry after A. Connes and A. L. Woronowicz. For
that reason we will use \cite{Connes:1986} and \cite{Woronowicz:1989} as the basic
references on the differential calculus on a Hopf algebra.

In recent papers, see~\cite{Aschieri/Castellani:1993}, \cite{Dubois:1995},
\cite{Jurco:1991}, \cite{Schmudgen:1996} and \cite{Wess/Zumino:1990}, the work of
Woronowicz was developed, in quantum groups and Hopf algebras, in order to obtain a non
commutative differential geometry and the correspondent (bi)covariant differential
calculus. Our purpose in this paper is to retrieve the notion of vector fields in this
noncommutative framework in such a way that they generalize non commutative vector
fields above mentioned and the vector fields as they was defined in the classical
(commutative) differential geometry.

Since some work in that direction, as we said before, was realized by P. Aschieri and
P. Schupp in ~\cite{Aschieri/Schupp:1995}, we recall their work and sketch the
relationship between our definition and their previous developments. The
Aschieri--Schupp approach consists in defining invariant vector fields from tangent
spaces and after that to build vector fields as the bimodule generated (over the Hopf
algebra) by invariant vector fields. Our approach here to vector fields is completely
different; we define vector fields from a global point of view using Cartan pairs, a
notion introduced by A. Borowiec, see~\cite{Borowiec:1996}, and after that we show that
the Aschieri--Schupp definition agrees with this new global definition.

The content of this paper is divided in four sections. In section one we define First
Order Differential Calculus and the construction of vector fields as they were
introduced by Aschieri and Schupp. In section two we introduce part of the theory of
Cartan pairs. In section three we study invariant vector fields and extend them to
covariant and non covariant First Order Differential Calculus. And the last section
deals with the definition of the Lie bracket of two invariant vector fields and their
properties.

Throughout this paper we use, as main reference for non defined notions, the
book~\cite{Klimyk/Schmudgen:1997} of A. Klimyk and K. Schm{\"u}dgen.

\section{Background. First Order Differential Calculus}

Let $R$ be an algebra. A \emph{First Order Differential Calculus}, FODC in brief, over $R$ is a
pair $(d,\Omega)$ where:
\newline(1)
$\Omega$ is an $R$--bimodule,
\newline(2)
$d\colon{R}\longrightarrow\Omega$ is a derivation, i.e., $d(ab)=d(a)b+ad(b)$ for any $a$,
$b\in{R}$, and
\newline(3)
$\Omega$ is generated, as $R$--bimodule, by the set $\{d(a)\colon\;a\in{R}\}$.

We may deduce, easily from the definition, that $\Omega$ is also generated, as left or
right $R$--module, by the set $\{d(a)\colon\;a\in{R}\}$ as for any $a$, $b\in{R}$ we
have the identity $d(a)b=d(ab)-ad(b)$.

\begin{example}
If $M$ is a $C^\infty$--differential manifold, we represent by $R=\mathcal{F}(M)$, the
commutative ring of either real or complex valued functions, and by
$\Omega=\Omega^1(M)$, the set of all 1-forms, then the pair $(d,\Omega)$ is a first
order differential calculus over $R$.
\end{example}

\begin{example}
Let $R$ be a $k$--algebra with multiplication map
$\mu\colon{R\otimes{R}}\longrightarrow{R}$. Then in the short exact sequence
\[
 \Omega^1_u(R)=\Ker(\mu)\hookrightarrow{R\otimes{R}}\stackrel{\mu}{\longrightarrow}R.
\]
$\Omega^1_u(R)$ is an $R$--bimodule as $\mu$ is an $R$--bimodule map.
\par
The linear map $d\colon{R}\longrightarrow\Omega^1_u(R)$, defined by
$d(a)=a\otimes1-1\otimes{a}$, is surjective. Indeed, let
$\sum{a_i\otimes{b_j}}\in\Ker(\mu)$, then we have
\[
\begin{array}{ll}
 \sum_{i}a_i\otimes{b_i}
 &=\sum_i(a_i\otimes1)b_i\\
 &=\sum_i(a_i\otimes1)b_i-\sum_i(1\otimes{a_i})b_i\\
 &=\sum_i(a_i\otimes1-1\otimes{a_i})b_i
\end{array}
\]
Thus $(d,\Omega^1_u(R))$ is a first order differential calculus over $R$.
\end{example}

\begin{lemma}
The pair $(d,\Omega^1_u(R))$ is a universal derivation, i.e., if $\Omega$ is an
$R$--bimodule and $D\colon{R}\longrightarrow\Omega$ is any derivation, then there
exists a unique linear map $f_D$ such that $D=f_D\circ d$.
\end{lemma}

As a consequence of this Lemma we call $(d,\Omega^1_u(R))$ the \emph{universal First
Order Differential Calculus} over $R$.

\begin{proposition}[{\cite{Woronowicz:1989}}]
Let $(d,\Omega^1_u(R))$ be the universal FODC over $R$. If $N$ is a subbimodule of
$\Omega^1_u(R)$ and we define
 $\Omega=\Omega^1_u(R)/N$,
 $\pi\colon\Omega^1_u(R)\rightarrow \Omega$ the canonical projection and
 $D=\pi\circ d$,
then $(D,\Omega)$ is a FODC over $R$.
\newline
In particular any FODC over $R$ can be obtained in this way.
\end{proposition}

If $(d,\Omega)$ is a FODC over $R$, elements in $\Omega$ are called \emph{1-forms}.

\subsection{First Order Differential Calculus over Hopf Algebras}

Following the classical differential geometric approach, the most interesting and
fruitful ground are Lie groups, i.e., manifolds with a differential group structure. In
this setting Hopf algebras appear as sets of real or complex valued functions over
compact Lie groups. Thus, from this and the above examples, Hopf algebras will be one
of the main structures in which the FODC must be studied. Let us recall the basic
definitions of FODC over Hopf algebras.

Let $H$ be a Hopf algebra, over a field $k$, and $(d,\Omega)$ be a FODC over $H$. We say that
$(d,\Omega)$ is \emph{left--covariant} if there exists a linear map
$\Delta_l\colon\Omega\rightarrow H\otimes\Omega$ such that the following diagram commutes.
\[
\begin{xy}
 \xymatrix{ H\ar[r]^d\ar[d]_\Delta & \Omega\ar@{-->}[d]_{\Delta_l}\ar@{=}[dr]\\
            H\otimes H\ar[r]_{I\otimes d} & H\otimes\Omega\ar[r]_{\varepsilon\otimes I} &\Omega}
\end{xy}
\]
From this commutativity it is easy to deduce that:
$\Delta_l(adb)=\Delta(a)(I\otimes{d})\Delta(b)$ for any $a$, $b\in H$.

In order to classify left--covariant FODC we need to associate certain elements. Let
$(d,\Omega)$ be a left--covariant FODC  over $H$. The set $\Omega_{inv}=\{\omega\in
\Omega\colon\Delta_l(\omega)=1\otimes\omega\}$ is the set of all \emph{invariant
1-forms}.
\par
We define $P_{inv}\colon\Omega\longrightarrow\Omega_{inv}$ by $P_{inv}(da)=S(a_1)da_2$,
where $S$ is the antipode in $H$ and $\Delta(a)=\sum a_1\otimes a_2$, in Sweedler's
notation and a right ideal $\mathcal{R}$ of $H$, contained in $\Ker(\varepsilon)$, as
follows
\[
\mathcal{R}=\{x\in\Ker(\varepsilon)\colon{P_{inv}(dx)=0}\}.
\]
From the Woronowicz's work the right ideals of $H$, contained in $\Ker(\varepsilon)$,
can be used to classify the left--covariant FODC. Let us show the Woronowicz's Theorem.

\begin{proposition}[{\cite{Woronowicz:1989}}]
Let $R$ be a right ideal of $H$ contained in $\Ker(\varepsilon)$ and
$N=\{aS(b_1)\otimes{b_2}\colon\;b\in{R},a\in H\}$ then $N$ is a subbimodule  of
$\Omega^1_u(H)$. Moreover, let $\Omega=\Omega^1_u(H)/N$ and let
$\pi\colon\Omega^1_u(H)\rightarrow \Omega$ be the canonical projection. if we define
$D=\pi\circ d$, then the FODC $(D,\Omega)$ over $H$ is left--covariant.
\newline
In particular any left--covariant FODC over $H$ can be obtained in this way.
\end{proposition}

\subsection{Vector fields}

We shall proceed to define the vector fields associated to a FODC. First remember, from
Lie group theory, that there exists a one-to-one correspondence between tangent vectors
on the neutral element $e$ of the Lie group, and invariant vector fields. Hence the
strategy, we could use to define vector fields, follows one of these parallel lines:
\newline\textbf{(i)}
 to use the tangent space on $e$ to define invariant vector fields or
\newline\textbf{(ii)}
 to define directly vector fields as ``\emph{derivations}" of $H$.

Definition of vector fields using (i) was developed by Schm{\"u}dgen--Sch{\"u}ler and
Aschieri--Schupp. We shall sketch in the following their approach.

Let $H'=\Hom(H,k)$ be the dual space of $H$ and let $(d,\Omega)$ be a left--covariant FODC over
$H$. The linear vector space
\[
 \chi=\{t\in H'\colon t(1)=0~\mbox{y}~t(x)=0~\forall~x\in\mathcal R\}
\]
is called the \emph{quantum (enveloping) Lie algebra} associated to the left--covariant FODC
$(d,\Omega)$. See~\cite[p. 3]{Schmudgen:1996}.

We recall that $H'=\Hom(H,k)$ is not in general a Hopf algebra, but the Hopf algebra
$H^o=\{f\in H'\colon f(I)=0,\mbox{ for some cofinite ideal }I\mbox{ of }H\}$ contains
the necessary information. In particular the counity $\varepsilon$ of $H$, is the unity
element of $H^o$.

\begin{lemma}[{\cite{Schmudgen:1996}}]
Let $\chi$ be a vector space of linear functionals on $H$, then the following statements are
equivalent:
\begin{enumerate}
\renewcommand{\theenumi}{(\alph{enumi})}
\item
$\chi$ is the quantum (enveloping) Lie algebra associated to a left--covariant FODC $(d,\Omega)$;
\item
$t(1)=0$ and $\Delta t-\varepsilon\otimes t\in\chi\otimes H'$ for all $t\in\chi$.
\end{enumerate}
\end{lemma}

In addition the FODC $(d,\Omega)$ may be obtained from the vector space $\chi$ as
follows. Let $\{t_i\colon i\in I\}$ be a basis of the vector space $\chi$ and
$\{x_i\colon i\in I\}$ be a set of elements of $H$ such that $t_i(x_j)=\delta_{ij}$. We
consider the invariant 1-forms $\omega_i=P_{inv}(dx_i)$ and define:
\[
da=\sum(t_i\ast a)\omega_i,~~a\in H
\]
where
\begin{equation}\label{smudprod}
f\ast a=(\mu\circ(I\otimes f)\circ\Delta)(a),
\end{equation}
see~\cite[p. 638]{Schmudgen:1995}. This product is denoted by $f\rightharpoonup a$
in~\cite{Montgomery:1993}

We may look at $\chi$ as the ``\emph{tangent space}" on $\varepsilon$. Next we shall
deduce from $\chi$ a definition of invariant vector fields. This was developed by
Aschieri and Schupp in \cite{Aschieri/Schupp:1995} as follows:
\[
\mathcal{X}_{inv}=\{V\in\Hom_k(H,H)\colon V=t\ast -~~ \mbox{ with
}t\in\chi\}
\]
The elements in $\mathcal{X}_{inv}$ are called \emph{invariant vector fields}. An element
$V\in\mathcal{X}_{inv}$ work as follows:
\[
V(a)=t\ast a=(I\otimes t)\circ\Delta (a)=a_1t(a_2)
\]

Also we can easily obtain $t$ from $V$ by composing with $\varepsilon$.

\begin{equation}\label{punto}
(\varepsilon\circ V)(a)=\varepsilon(a_1t(a_2))=t(a)
\end{equation}

Now, using (\ref{punto}), we may rewrite the definition of invariant vector field
without using elements in $\chi$.
\begin{proposition}
Let $V\in\Hom_k(H,H)$ then the following statements are
equivalent:
\begin{enumerate}
\renewcommand{\theenumi}{(\alph{enumi})}
\item $V$ is an invariant vector field,
\item $V$ satisfies:
\begin{equation}\label{defA}
 V=(I\otimes\varepsilon)(I\otimes V)\Delta
\end{equation}
\end{enumerate}
\end{proposition}

From a geometrical point of view we may look the composition with $\varepsilon$ as
calculate the value of vector field on the point $\varepsilon$, the unity of $H^o$. So,
formula (\ref{punto}) say us that if we construct an invariant vector field
$V=t\ast(-)$, from a ``vector" $t$ in the tangent space on $\varepsilon$, we can
recover $t$ in computing the value of $V$ at the point $\varepsilon$.

Let $\{t_i\colon i\in I\}$ be a basis of $\chi$, the quantum Lie algebra associated to
a left--covariant FODC $(d,\Omega)$. A general \emph{vector field} was defined as:
\[
X=\sum a^iV_i
\]
where $a^i\in H$ and $V_i=t_i\ast -~\in\mathcal{X}_{inv}$.

The disadvantage of this approach coming from the following fact: ``\emph{There are
some problems to perform right multiplication in this definition}". It is very
difficult to work with these vector fields as we can't controle the Leibnitz rule.
Remember that for $X=\sum a^iV_i$ Leibnitz rule says that
\[
X(ab)=X(a)b+a^i\sum_j(f_{ij}\ast a)t_j(b),
\]
being $\{t_i\}$ a basis of $\chi$ with $V_i$ the invariant vector fields associated,
and:
\[
t_i(ab)=t_i(a)\varepsilon(b)+\sum_j f_{ij}(a)t_j(b)
\]
being $f_{ij}$ functionals over $H$. See~\cite[Theorem 2.1]{Woronowicz:1989} for
details. We propose a new and natural approach to define vector fields and show that
they satisfies a Leibnitz rule in a more easy way. After that, the theory could be
extended to consider connections and other tools in the framework of non commutative
differential geometry.

\section{Cartan Pairs}

As a consequence of the difficulties expounded in the above section, we shall explore a
new way to define vector fields. We shall produce a slight modification: we consider
``\emph{derivations}" of the algebra $R$ as elements in $\Hom(R,R)$, instead of in
$R'=\Hom(R,k)$, following the second alternative announced in Section 1.

To do that we need the notion of Cartan pair as it was introduced by Borowiec in
\cite{Borowiec:1996}. Let $R$ be a $k$--algebra. A \emph{right Cartan pair} is a pair
$(M,\rho)$ satisfying the following properties:
 \bara
 \item $M$ is a bimodule;
 \item $\rho\colon{M}\longrightarrow\End_k(R)$ is a linear map;
 \item $\rho(a\cdot{m})(r)=a\cdot\rho(m)(r)$ ($\rho$ is a left $R$--module map);
 \item $\rho(m)(rs)=\rho(m)(r)s+\rho(mr)(s)$ ($\rho$ is a ``new derivation");
 \eara

A right Cartan pair is \emph{faithful} or a \emph{right vector bundle} if $\rho$ is an injective
map.

Let $(M_1,\rho_1)$ and $(M_2,\rho_2)$ be two right Cartan pairs, a
\emph{right Cartan pair map} from $(M_1,\rho_1)$ to $(M_2,\rho_2)$
is a linear map $g\colon{M_1}\longrightarrow{M_2}$ satisfying
 \bara
 \item The diagram
 \[
 \begin{xy}\xymatrix{
  M_1\ar[rr]^g\ar[rd]_{\rho_1}&&M_2\ar[ld]^{\rho_2}\\
  &\End_k(R)\\
 }\end{xy}
 \]
 is commutative and
 \item $g$ is a bimodule map.
 \eara

The relationship between Cartan pairs and FODC is given in the next theorem.

\begin{theorem}
Each FODC determines a faithful right Cartan pair.
\end{theorem}
\begin{proof}
Let $(d,\Omega)$ be a FODC. We consider $\Omega^\ast=\Hom_{-R}(\Omega,R)$, the set of all
right $R$--module homomorphisms and define
$\rho_d\colon\Omega^\ast\longrightarrow\End_k(R)$ as $\rho_d(X)(r)=X(dr)$.

First we obtain that $\Omega^\ast$ is a bimodule as for any $a$,
$b\in{R}$ and $\omega\in\Omega$ we have:
\[
(aXb)(\omega)=aX(b\omega).
\]
The map $\rho_d$ is linear. It satisfies property (3), i.e.,
$\rho_d(aX)(r)=(aX)(dr)=aX(dr)=a\rho_d(X)(r)$; and  property (4)
i.e.,
\[
\begin{array}{ll}
 \rho_d(X)(rs)
 &=X(d(rs))=X((dr)s+r(ds))\\
 &=X((dr)s)+X(r(ds))=X(dr)s+(Xr)(ds)\\
 &=\rho_d(X)(r)s+\rho_d(Xr)(s)
\end{array}
\]
as $X$ is a right $R$--module map. In addition if $\rho_d(X)=0$, then $\rho_d(X)(s)=0$ for any
$s\in{R}$, i.e, $X(ds)=0$ for any $s\in{R}$. If $d(R)$ spans $\Omega$ as bimodule, then we have:
\[
\begin{array}{ll}
 X(a(dr)b)
 &=X(a(dr))b\\
 &=X(d(ar)-(da)r)\\
 &=X(d(ar))b-X(da)rb\\
 &=\rho_d(X)(ar)b-\rho_d(X)(a)rb=0.
\end{array}
\]
Therefore $X$ is zero.
\end{proof}

\section{Translations on Hopf Algebras}

In order to compare both definitions, we shall use invariant vector fields. We justify
this choice in the following: the invariant vector fields constitute a set of
generators for the $R$-bimodule of all vector fields, hence we only need to know these
invariant vector fields to obtain all vector fields, and remember that it is much
easier to work with invariant vector fields than with vector fields.

In order to establish a definition of invariance for vector fields, we translate to
Hopf algebras some ideas from the classical differential geometry.

The invariance notion in differential geometry is associated with left translations,
hence we define translation in a Hopf algebra framework as follows.

Let $H$ be a Hopf algebra. We say that $L\colon H\rightarrow H$ is a \emph{left
translation} on $H$, if $L$ is an algebra map and the following diagram commute:
\[
\begin{xy}
\xymatrix{H\ar[r]^{\Delta}\ar[d]_L & H\otimes H\ar[d]^{L\otimes I}\\
  H\ar[r]_{\Delta} & H\otimes H}
\end{xy}
\]
i.e., $L$ is a right comodule map where the right comodule structure on $H$ is given
via $\Delta$.

Right translations can be defined as algebra maps from $H$ to $H$ being left comodule
maps. In the following we will use simply translation instead of left translation.

It is clear that the composition of translations is again a translation.

\begin{proposition}\label{trascom}
Let $\lambda\colon C\rightarrow H\otimes C$ be a comodule structure on $H$, i.e.,
$(\Delta\otimes I)\lambda=(I\otimes\lambda)\lambda$, and $L\colon H\rightarrow H$ be a
translation, then we can define a linear map $L'\colon C\rightarrow C$ satisfying that
the following diagram commutes:
\[
\begin{xy}
\xymatrix{C\ar[r]^{\lambda}\ar[d]_{L'} & H\otimes C\ar[d]^{L\otimes I}\\
  C\ar[r]_{\lambda} & H\otimes C}
\end{xy}
\]
\end{proposition}
\begin{proof}
In the following the diagram the two compositions
$(1\otimes\lambda)(L\otimes{I})\lambda$ and $(\Delta\otimes{I})(L\otimes{I})\lambda$
are equal.
\[
\begin{xy}
\xymatrix{ C\ar[r]^{\lambda} &  H\otimes C\ar[d]^{L\otimes
I}\ar[r]^{\Delta\otimes I} &
                              H\otimes H\otimes C\ar[d]^{L\otimes I\otimes I}\\
   C\ar[r]_{\lambda} &  H\otimes  C\ar[r]<3pt>^{\Delta\otimes I}\ar[r]<-3pt>_{I\otimes\lambda} &
                              H\otimes  H\otimes C}
\end{xy}
\]
Indeed, writing $\lambda(c)=c_1\otimes c_2$ in Sweedler notation's, we have:
\[
\begin{array}{rcl}
 ((I\otimes\lambda)\circ(L\otimes I)\circ \lambda) (c)
                & = & ((I\otimes\lambda)\circ(L\otimes I))(c_1\otimes c_2)\\
                & = & (I\otimes\lambda)(Lc_1\otimes c_2) \\
                & = & Lc_1\otimes c_2\otimes c_3 \\
((\Delta\otimes I)\circ(L\otimes I)\circ \lambda) (c)
                & = & ((\Delta\otimes I)\circ(L\otimes I))(c_1\otimes c_2)\\
                & = & (\Delta\otimes I)(Lc_1\otimes c_2) \\
                & = & (\Delta\circ L) c_1\otimes c_2 \\
                & = & ((L\otimes I)\circ\Delta)(c_1)\otimes c_2 \\
                & = & Lc_1\otimes c_2\otimes c_3
\end{array}
\]
For any $c\in{C}$ we define a new element $\alpha:=(\varepsilon\otimes I)\circ(L\otimes
I)\circ\lambda(c)\in{C}$ and show that it satisfies $\lambda(\alpha)=(L\otimes
I)\circ\lambda(c)$. Then the map $L'\colon{C}\to{C}$, defined by $L'=(\varepsilon\otimes
I)\circ(L\otimes I)\circ\lambda(c)$, is an answer to the Proposition.
\end{proof}

We are now interested on the geometrical meaning of $\varepsilon\circ L$. First we
observe that $\varepsilon\circ L$ is an algebra map as both, $\varepsilon$ and $L$,
are.

\begin{proposition}
Let $L$ be a translation, then there exists an algebra map $\varphi\in\Alg(H,k)$ such
that $(\varphi\otimes I)\circ\Delta=L$.
\end{proposition}
\begin{proof}
Just define $\varphi=\varepsilon\circ L$.
\end{proof}

\begin{remark}
If $R$ is right translation, there exists $\varphi'\in\Alg(H,k)$
such that $(I\otimes\varphi')\circ\Delta=R$.
\end{remark}

We may prove that the converse of the above Proposition is also true.

\begin{proposition}
Let $\varphi\in\Hom(H,k)$ then $L=(\varphi\otimes I)\circ\Delta$ is a right comodule map.
Moreover, $L$ is an algebra map if and only if $\varphi$ is an algebra map.
\end{proposition}
\begin{proof}
In fact, we only need to perform the following computation:
\[
\begin{array}{rcl}
 \Delta\circ L & = & \Delta\circ(\varphi\otimes I)\circ\Delta \\
               & = & (\varphi\otimes I\otimes I)\circ(I\otimes\Delta)\circ\Delta\\
               & = & (\varphi\otimes I\otimes I)\circ(\Delta\otimes I)\circ\Delta\\
               & = & (L\otimes I)\circ\Delta
\end{array}
\]
\end{proof}

Our next goal is, using the last equivalence, to relate the structure of the two
following sets:
 (1) the set of all left translations $L:H\to{H}$ and
 (2)the structure of the set of algebra maps $\varphi:H\to{k}$.
To do that, first we remember that the algebra maps from $H$ to $k$ are parameterized
by the group--like elements of the Hopf algebra $H^o$. Hence we identify $\Alg(H,k)$
and $G(H^o)$.

In $G(H^o)$ we may define a ``\emph{new multiplication}", represented by $\star$, using
the composition of translations and the above identification. Thus if $\varphi$,
$\varphi'\in\Alg(H,k)=G(H^o)$, with associated translation $L$ and $L'$ respectively,
then
\[
 \varphi\star\varphi'=\varepsilon(L\circ{L'}).
\]
The next Theorem related this new multiplication with the \emph{anticonvolution} or
\emph{twist convolution} multiplication. See~\cite{Montgomery:1993} for the
definitions.

\begin{theorem}
The multiplication $\star$ in $G(H^o)$, coincides with the \emph{anticonvolution} or
\emph{twist convolution} multiplication, i.e.,
\[
\varphi\star\varphi'=\varepsilon(L\circ
L')=\mu\circ(\varphi\otimes\varphi')\circ\tau\circ\Delta,
\]
where $\tau$ is the \emph{twist map} ($\tau(a\otimes b)=b\otimes a$).
\newline
In this case, if $S$ is the antipode, then the inverse of $\varphi$ is the composition
$\varphi\circ{S}$.
\end{theorem}
\begin{proof}
We write $L=(\varphi\otimes I)\circ\Delta$ and $L'=(\varphi'\otimes I)\circ\Delta$.
Then we obtain:
\[
\begin{array}{rcl}
 (L\circ L')(a) & = & ((\varphi\otimes I)\circ\Delta\circ (\varphi'\otimes I)\circ\Delta)(a)\\
                & = & ((\varphi\otimes I)\circ\Delta)(\varphi'(a_1)a_2)\\
                & = & \varphi'(a_1)\varphi(a_2)a_3
\end{array}
\]
If we compose with $\varepsilon$ then
\[
\begin{array}{rcl}
 \varepsilon\circ(L\circ L')(a)
   & = & \varepsilon(\varphi(a_2)\varphi'(a_1)a_3)\\
   & = & \varphi(a_2)\varphi'(a_1)\varepsilon(a_3)\\
   & = & \varphi(a_2\varepsilon(a_3))\varphi'(a_1)\\
   & = & (\mu\circ(\varphi\otimes\varphi')\circ\tau\circ\Delta)(a)
\end{array}
\]
and the first assertion holds.

To finish the proof, we recall some basic facts on the antipode $S$ of $H$ and $H^o$.
\newline (i)
If $\varphi\in\Alg(H,k)$ then $\varphi\circ{S}\in\Alg(H,k)$ as $S:H\to{H}$ is an
anti--algebra homomorphism.
\newline (ii)
The antipode $S^*$ of $H^o$ is defined by $S^*(f)=f\circ{S}$ for any $f\in{H^o}$
(see~\cite[theorem 9.1.3]{Montgomery:1993}).
\newline (iii)
For any group--like element $\varphi$ of $H^o$ we have that $S^*$ works as
$S^*(\varphi)=\varphi^{-1}$ (see~\cite[example 1.5.3]{Montgomery:1993}).

We deduce from (ii) and (iii) that for any $\varphi\in\Alg(H,k)$ we obtain
\begin{equation}\label{inversa}
\varphi\circ S=\varphi^{-1}.
\end{equation}
Now, using (i) we have $\varphi^{-1}\circ S\in G(H^o)$ and as before we have
$\varphi^{-1}\circ S=\varphi$.

We may also compute directly the formula (\ref{inversa}) as follows:
\[
\begin{array}{rcl}
 (\varphi\circ S)\star \varphi(a) & = & (\varphi\circ S)(a_2)\cdot\varphi(a_1)\\
                                   & = & \varphi(a_1)\cdot\varphi(S(a_2))\\
                                   & = & \varphi(\varepsilon(a)\cdot 1) \\
                                   & = & \varepsilon(a)
\end{array}
\]
\end{proof}

The relationship between this product and the product defined in (\ref{smudprod}) is
given by the following formula:
\[
(\varphi\star\varphi')(a)=\varphi\circ(\varphi'\ast a).
\]
As it is obvious, the antipode $S$ relates the right and left sides of the Hopf
algebra, and many results can be established using $S$. Let us point out the following
one for future applications.

\begin{proposition}
If $L:H\to{H}$ is a right comodule map, then there exists a left comodule map
$R:H\to{H}$ such that $R\circ S=S\circ L$.
\end{proposition}
\begin{proof}
From $L$ we define a new map as follows: $R=(I\otimes(\varepsilon\circ
L^{-1}))\circ\Delta$. It is easy to check that $R$ is a left comodule map. In order to
check that $R$ and $L$ are related by the formula, we proceed as follows:
\[
\begin{array}{rcl}
 (R\circ S)(a) & = & ((I\otimes(\varepsilon\circ L^{-1}))\circ\Delta\circ S)(a) \\
             & = & S(a_2)\cdot(\varepsilon\circ L^{-1}\circ S)(a_1) \\
             & = & S(a_2)\cdot(\varphi^{-1}\circ S)(a_1) \\
             & = & S(a_2)\cdot\varphi(a_1)\\
             & = & S(\varphi(a_1)\cdot a_2)\\
             & = & S(L(a))=(S\circ L)(a)
\end{array}
\]
\end{proof}

\subsection{Translations over Covariant First Order Differential
Calculus}

Once we have established the basic behavior of translations, we shall apply them to
study covariant FODC on Hopf algebras. Thus the translation $L$ defines, in a natural
way, a new map $L'\colon\Omega\longrightarrow\Omega$.

\begin{proposition}
Let $(d,\Omega)$ be a covariant FODC over a Hopf algebra $H$. For every translation $L$
we can find $L'\colon\Omega\rightarrow\Omega$ such that the following diagram commutes.
\[
\begin{xy}
\xymatrix{\Omega\ar[r]^{\Delta_l}\ar[d]_{L'} & H\otimes \Omega\ar[d]^{L\otimes I}\\
  \Omega\ar[r]_{\Delta_l} & H\otimes \Omega}
\end{xy}
\]
\end{proposition}
\begin{proof}
We can obtain the result applying Proposition \ref{trascom} for
$\Delta_l\colon\Omega\rightarrow H\otimes\Omega$. So it's only necessary to check that
$(\Delta\otimes I)\Delta_l=(I\otimes\Delta_l)\Delta_l$.
\par
Indeed, we compute the following expression.
\[
\begin{array}{ll}
 ((\Delta\otimes I)\circ\Delta_l)(adb)
 &=(\Delta\otimes I)(a_1b_1\otimes a_2db_2)\\
 &=a_1b_1\otimes a_2b_2\otimes a_3db_3 \\
 &=(I\otimes\Delta_l)(a_1b_1\otimes a_2db_2)\\
 &=((I\otimes\Delta_l)\circ\Delta_l) (adb)
\end{array}
\]
\end{proof}

As a matter of fact, we may define $L'=(\varepsilon\otimes I)\circ(L\otimes
I)\circ\Delta_l$. Hence we obtain the following relationship for $L'$.
\[
\begin{array}{rcl}
 &&L'(a\omega b)\\
 & = & ((\varepsilon\otimes I)\circ(L\otimes I)\circ\Delta_l)(a\omega b)\\
 & = & ((\varepsilon\otimes I)\circ(L\otimes I))(\Delta(a)\Delta_l(\omega)\Delta(b))\\
 & = & ((\varepsilon\otimes I)\circ\Delta\circ L)(a)\cdot((\varepsilon\otimes I)\circ(L\otimes I)\circ\Delta_l)(\omega)\cdot((\varepsilon\otimes I)\circ\Delta\circ L)(b)\\
 & = & L(a)L'(\omega)L(b)
\end{array}
\]
But it is possible to obtain an easy definition for $L'$. Indeed, we may build the
following diagram:
\[
\begin{xy}
 \xymatrix{ & H\ar@{-->}[rr]^d \ar@{-->}[dd]^(0.4)L\ar[dl]_\Delta & & \Omega\ar[dr]^{\Delta_l}\ar@{-->}[dd]^(0.4){L'}  & \\
   H\otimes H \ar[dd]_{L\otimes I}\ar[rrrr]^{I\otimes d} & & & & H\otimes\Omega \ar[dd]^{L\otimes I} \\
            & H\ar@{-->}[rr]^d \ar[dl]_\Delta &  & \Omega\ar[dr]^{\Delta_l}  & \\
   H\otimes H \ar[rrrr]^{I\otimes d} & & & & H\otimes\Omega}
\end{xy}
\]
Where all the squares, except possibly, the square with broken arrows, are commutative.
But also this square is commutative.

\begin{theorem}
In the above diagram the square with broken arrows is commutative, i.e., $L'd=dL$.
\end{theorem}
\begin{proof}
Indeed, we have:
\[
\begin{array}{ll}
 \Delta_l\circ L'\circ d
 &=(L\otimes I)\circ\Delta_l\circ d\\
 &=(L\otimes I)\circ (I\otimes d)\circ\Delta\\
 &=(I\otimes d)\circ(L\otimes I)\circ \Delta\\
 &=(I\otimes d)\circ\Delta\circ L=\Delta_l\circ d\circ L
\end{array}
\]
then, after composing with $(\varepsilon\otimes I)$, we finish the proof.
\end{proof}

\subsection{Invariant Vector Fields Using Translations}

We define a \emph{vector field} as an element $X$ in $\Omega^*=\Hom_{-,H}(\Omega,H)$,
where $(d,\Omega)$ is a FODC over $H$.

For any translation $L$ we have defined a map $L'\colon\Omega\rightarrow \Omega$, and
this map induces a map $L^*\colon\Omega^*\rightarrow\Omega^*$ by $(L^*
X)(\omega)=X(L'\omega)$, which is related with $\rho(X)$ by the following relationship:
\[
\begin{array}{ll}
 \rho(L^*X)(a)
 &=(L^*X)(da)\\
 &=X(L'da)\\
 &=X(dLa)\\
 &=\rho(X)(La)\\
 &=(\rho(X)\circ L) (a)
\end{array}
\]

A vector field $X\in\Omega^*$ is an \emph{invariant vector field} if $L^*(X)=L\circ X$
for any translation $L$.

Using $\rho$ we may also write
\begin{equation}\label{defB}
\rho(X)\circ L=L\circ \rho(X),
\end{equation}
as $\rho(X)\circ L=\rho(L^*X)$ and
\[
\rho(L\circ X)(a)=(L\circ X)(da)=L(X(da))=L\circ\rho(X)(a).
\]
Relation (\ref{defB}) may be written as:
\[
L\circ \rho(X)=\rho(X)\circ(\varphi\otimes I)\circ\Delta=(\varphi\otimes
I)\circ(I\otimes\rho(X))\circ\Delta,
\]
and if we compose with $\varepsilon$ then for all $\varphi\in\Alg(H,k)$ we obtain
\[
\varphi\circ\rho(X)=\varepsilon\circ(\varphi\otimes
I)\circ(I\otimes\rho(X))\circ\Delta=\varphi\circ(I\otimes
\varepsilon)\circ(I\otimes\rho(X))\circ\Delta.
\]

Let us consider a particular case. If we assume that $\Alg(H,k)$ separates the elements
in $H$, i.e., given $a$, $b\in H$, $a\neq b$, exists $\varphi\in\Alg(H,k)$ such that
$\varphi(a)\neq\varphi(b)$, then condition (\ref{defB}) implies condition (\ref{defA}).

The converse is also true. Indeed, as $\rho(X)=(I\otimes \varepsilon)\circ(I\otimes
\rho(X))\circ\Delta$, after composing with $L$ we have:
\[
\begin{array}{rcl}
  \rho(X)\circ L & = & (I\otimes \varepsilon)\circ(I\otimes \rho(X))\circ\Delta\circ L \\
                 & = & (I\otimes \varepsilon)\circ(I\otimes \rho(X))\circ(L\otimes I)\circ\Delta \\
                 & = & (I\otimes \varepsilon)\circ(L\otimes I)\circ(I\otimes \rho(X))\circ\Delta \\
                 & = & L\circ(I\otimes \varepsilon)\circ(I\otimes \rho(X))\circ\Delta  \\
                 & = & L\circ\rho(X)
\end{array}
\]
So, in that case, definitions (\ref{defB}) and (\ref{defA}) are equivalent.

\section{Lie Bracket for Invariant Fields}

Recall that in commutative differential geometry, the Lie bracket defines the structure
of Lie algebra of the set of vector fields, and it is well known how this construction
is used in the study of differential varieties. The problem that arise is how to
realize this construction in non commutative differential geometry. The first approach
to the construction of the Lie bracket in a noncommutative case was realized
in~\cite{Dubois:1995} by considering elements in the center of $R$ (in that case all
the elements commute). Our purpose in this section is to extend these results and
define a new Lie bracket in such a way that it works for invariant vector fields.

We start this section looking for vector fields satisfying the
property required in the definition of Cartan pair.

First we realize, in the next paragraph, a change of language.

Let $(M,\rho)$ be a Cartan pair over an algebra $R$, let $X$ be a vector field and
$a\in{R}$. If we write $\rho(Xa)=\sum\rho_1(a)\cdot\rho_2$, then we obtain
$\Delta(\rho(X))=\sum\rho_1\otimes\rho_2+\rho(X)\otimes I$. Indeed, we have:
\[
 \Delta(\rho(X))(a\otimes b)
 =\rho(X)(ab)
 =(\sum\rho_1\otimes\rho_2+\rho(X)\otimes I)(a\otimes b)
\]

Let us analyze a particular case which is of importance in the sequel.

\begin{lemma}
If the following identity holds
\begin{equation}\label{covdif}
\Delta(\rho(X))=L\otimes\rho(X)+\rho(X)\otimes I,
\end{equation}
then $Xa=L(a)X$.
\end{lemma}
\begin{proof}
Since $\rho(X)(ab)=\rho(Xa)(b)+\rho(X)(a)\cdot b$; and using the
hypothesis we deduce:
$\rho(X)(ab)=L(a)\cdot\rho(X)(b)+\rho(X)(a)\cdot b$. So
$Xa(db)=L(a)\cdot X(db)$.
\end{proof}

The election of vector fields satisfying property (\ref{covdif}) should be no strange,
as in the commutative case the vector fields satisfy
$\Delta(\rho(X))=I\otimes\rho(X)+\rho(X)\otimes I$.

\begin{theorem}
Let $X$ and $Y$ be two invariant vector fields, i.e., they commute with all
$L\in\Alg(R)$, satisfying:
\[
\begin{array}{rcl}
\Delta(\rho(X)) & = & L\otimes\rho(X)+\rho(X)\otimes I\quad\mbox{ and} \\
\Delta(\rho(Y)) & = & L'\otimes\rho(Y)+\rho(Y)\otimes I.
\end{array}
\]
If we assume
\begin{equation}\label{comm}
L\circ L'=L'\circ L,
\end{equation}
then we have:
\[
\rho([X,Y])=\rho(X)\circ\rho(Y)-\rho(Y)\circ\rho(X).
\]
Moreover $[X,Y]$ is an invariant vector field, i.e.,
\[
 \Delta(\rho([X,Y]))=L\circ L'\otimes\rho([X,Y])+\rho([X,Y])\otimes I.
\]
\end{theorem}
\begin{proof}
In fact, after applying $\Delta$ we obtain:
\[
\begin{array}{rcl}
 \Delta(\rho([X,Y])) & = & \Delta(\rho(X)\circ\rho(Y)-\rho(Y)\circ\rho(X)) \\
                     & = & L\circ L'\otimes\rho(X)\circ\rho(Y) +\rho(X)\circ L'\otimes\rho(Y)\\
                     &  & +L\circ\rho(Y)\otimes\rho(X)+\rho(X)\circ\rho(Y)\otimes I\\
                     &  & -L'\circ L\otimes\rho(Y)\circ\rho(X) -\rho(Y)\circ L\otimes\rho(X)\\
                     &  & -L'\circ\rho(X)\otimes\rho(Y) -\rho(Y)\circ\rho(X)\otimes I\\
                     & = & L\circ L'\otimes\rho([X,Y])+\rho([X,Y])\otimes I
\end{array}
\]
The invariance property is easy to check.
\end{proof}

It's important, in this proof, the invariance of the vector fields $X$ and $Y$, because
otherwise the definition of bracket don't work see~(\ref{noinvariant}) below.

Invariant vector fields satisfying property (\ref{covdif}) will be called
\emph{invariant coderivations}.

We point out that if in the above Theorem we don't use $\rho$, then we may write:
\[
[X,Y](da)=X(dY(da))-Y(dX(da))
\]
as in the classical (commutative) differential geometry.

\begin{proposition}
The Leibnitz rule for invariant coderivations satisfying (\ref{comm}) is wrote, using
the last notation, as follows:
\[
[X,Y](d(ab))  =  L\circ L'(a)\cdot [X,Y](db)+[X,Y](da)\cdot b
\]
\end{proposition}
\begin{proof}
By definition we have
\[
[X,Y](dab)=X(dY(dab))-Y(dX(dab)),
\]
hence:
\[
\begin{array}{lcl}
 X(dY(dab))
 &=&X(dY(adb+da\cdot b))\\
 &=&X(dY(adb))+X(dY(da\cdot b))\\
 &=&X(d(L'(a)\cdot Y(db)))+X(d(Yda)\cdot b))\\
 &=&X(L'(a)\cdot d(Y(db)))+X(d(L'(a))\cdot Y(db))\\
 & & + X(Yda\cdot db))+X(d(Yda)))\cdot b)\\
 &=& L\circ L'(a)X(d(Y(db)))+L'(X(da))\cdot Y(db)\\
 & & +L(Yda)X(db)+X(d(Yda))\cdot b.
\end{array}
\]
In an analogous way we have:
\[
\begin{array}{rcl}
 Y(dX(dab)) & = & L'\circ L(a)Y(d(X(db)))+L(Y(da))\cdot X(db)\\
            &   & +L'(Xda)Y(db)+Y(d(Xda))\cdot b.
\end{array}
\]
If we reduce both, using the relationship $L\circ L'=L'\circ L$, then we have the
result.
\end{proof}

\subsection{Others Computations with the Bracket}

Let us give two more formules to see the importance of the invariance property for
coderivations. Recall that if $X$ is an invariant vector field, then $Xa$ is invariant
only if $a\in k$.

\begin{theorem}
Let $X$, $Y$ be two invariant coderivations satisfying
(\ref{comm}) and be $a$, $b\in R$ then we have:
\begin{itemize}
 \item
 $
 [X,Y](adb)  =  L\circ L'(a)[X,Y](db)+L'(Xda)\cdot Y(db)-L(Yda)\cdot X(db)
 $
 \item
 $
 [X,Y](da\cdot b)  = [X,Y](da)\cdot b + L(Yda)\cdot X(db)-L'(Xda)\cdot Y(db)
 $
\end{itemize}
\end{theorem}
\begin{proof}
The first assertion follows from the following computations:
\[
\begin{array}{rcl}
 X(dY(adb)) & = & X(d(L'(a)\cdot Ydb))\\
            & = & X(L'(a)\cdot d(Ydb))+ X(d(L'(a))\cdot Ydb)\\
            & = & L\circ L'(a)\cdot X(d(Ydb))+L'(Xda)\cdot Y(db).
\end{array}
\]
On the other hand we have:
\[
 Y(dX(adb))= L'\circ L(a)\cdot Y(d(Xdb))+L(Yda)\cdot X(db),
\]
the result holds.

To see the second assertion we compute in the same way and obtain:
\[
\begin{array}{rcl}
 X(d(Yda\cdot b))  & = & X(Yda\cdot db))+X(d(Yda)\cdot b)\\
                   & = & L(Yda)\cdot X(db)+ X(d(Yda))\cdot b
\end{array}
\]
and
\[
 Y(d(Xda\cdot b))= L'(Xda)\cdot Y(db)+Y(d(Xda))\cdot b.
\]
Thus we have the result.
\end{proof}

As a final result let us enumerate, in the following Theorem, the properties satisfied
by the Lie bracket.

\begin{theorem}\label{noinvariant}
We have too for $X,Y,Z$ invariant coderivations satisfying
(\ref{comm}) and $a,b,c\in R$
\begin{itemize}
 \item $[X,Y]=-[Y,X]$
 \item $[X+X',Y]=[X,Y]+[X',Y]$
 \item $[X,[Y,Z]]+[Z,[X,Y]]+[Y,[Z,X]]=0$
 \item $[aX,Y](db)=aX(dY(db))-L'(a)Y(d(Xdb))-Y(da)\cdot X(db)$
 \item $\begin{array}{lcl}
        [Xa,Y](db)
        &=&L(a)X(dY(db))-L'\circ L(a)Y(d(Xdb))\\
        & &-L(Y(da))\cdot X(db)
        \end{array}$
 \item $[X,aY](db)=L(a)X(d(Ydb))+X(da)\cdot Y(db)-a\cdot Y(d(Xdb))$
 \item $\begin{array}{lcl}
        [X,Ya](db)
        &=&L\circ L'(a)X(d(Ydb))+L'(Xda)\cdot Y(db)\\
        & &-L'(a)\cdot Y(d(Xdb))
        \end{array}$
 \item $\begin{array}{lcl}
        [aX,bY](dc)
        &=&aL(b)X(d(Ydc))+a(Xdb)(Ydc)\\
        & &-bL'(a)Y(d(Xdc))-b(Yda)(Xdc)
        \end{array}$
\end{itemize}
\end{theorem}

All these formulas are easy to check and we left them as an exercise to the reader.

\end{document}